\documentclass[11pt]{amsart}
\usepackage{amsmath,amssymb}
\newtheorem{theorem}{Theorem}
\newtheorem{proposition}[theorem]{Proposition}
\newtheorem{corollary}[theorem]{Corollary}
\newtheorem{lemma}[theorem]{Lemma}
\theoremstyle{definition}

\begin{document}

\title[Michael-Simon-Sobolev inequality]{Proof of the Michael-Simon-Sobolev inequality using optimal transport}
\author{Simon Brendle and Michael Eichmair}
\address{Columbia University, 2990 Broadway, New York NY 10027, USA}
\address{University of Vienna, Faculty of Mathematics, Oskar-Morgenstern-Platz 1, 1090 Vienna, Austria, ORCID: 0000-0001-7993-9536}
\begin{abstract}
We give an alternative proof of the Michael-Simon-Sobolev inequality using techniques from optimal transport. The inequality is sharp for submanifolds of codimension $2$.
\end{abstract}
\thanks{The first-named author was supported by the National Science Foundation under grant DMS-2103573 and by the Simons Foundation. The second-named author was supported by the START-Project Y963 of the Austrian Science Fund.}

\maketitle

\section{Introduction} 

In this paper, we use techniques from optimal transport to prove the following result.

\begin{theorem} 
\label{main.theorem}
Let $n \geq 2$ and $m \geq 1$ be integers. Let $\rho: [0,\infty) \to (0,\infty)$ be a continuous function with $\int_{\bar{B}^{n+m}} \rho(|\xi|^2) \, d\xi = 1$, where $\bar{B}^{m+m} = \{\xi \in \mathbb{R}^{n+m}: |\xi| \leq 1\}$ denotes the closed unit ball in $\mathbb{R}^{n+m}$. Let 
\begin{equation} 
\label{def.alpha} 
\alpha = \sup_{z \in \mathbb{R}^n} \int_{\{y \in \mathbb{R}^m: \, |z|^2+|y|^2 \leq 1\}} \rho(|z|^2+|y|^2) \, dy. 
\end{equation}
Let $\Sigma$ be a compact $n$-dimensional submanifold of $\mathbb{R}^{n+m}$, possibly with boundary $\partial \Sigma$. Then 
\begin{equation} 
\label{michael.simon}
|\partial \Sigma| + \int_\Sigma |H| \geq n \, \alpha^{-\frac{1}{n}} \, |\Sigma|^{\frac{n-1}{n}}, 
\end{equation}
where $H$ denotes the mean curvature vector of $\Sigma$.
\end{theorem} 

The proof of Theorem \ref{main.theorem} is based on an optimal mass transport problem between the submanifold $\Sigma$ and the unit ball in $\mathbb{R}^{n+m}$, the latter equipped with a rotationally invariant measure. A notable feature is that this transport problem is between spaces of different dimensions. 

In Theorem \ref{main.theorem}, we are free to choose the density $\rho$. For $m \geq 2$, it is convenient to choose the density $\rho$ so that nearly all of the mass of the measure $\rho(|\xi|^2) \, d\xi$ on $\bar{B}^{n+m}$ is concentrated near the boundary. This recovers the main result of \cite{Brendle1}.

\begin{corollary} 
\label{consequence.of.main.theorem}
Let $n \geq 2$ and $m \geq 2$ be integers. Let $\Sigma$ be a compact $n$-dimensional submanifold of $\mathbb{R}^{n+m}$, possibly with boundary $\partial \Sigma$. Then 
\begin{equation} 
\label{michael.simon.2}
|\partial \Sigma| + \int_\Sigma |H| \geq n \, \Big ( \frac{(n+m) \, |B^{n+m}|}{m \, |B^m|} \Big )^{\frac{1}{n}} \, |\Sigma|^{\frac{n-1}{n}},
\end{equation}
where $H$ denotes the mean curvature vector of $\Sigma$.
\end{corollary} 

Note that the constant in (\ref{michael.simon.2}) is sharp for $m=2$. 

Earlier proofs of the non-sharp version of the inequality were obtained by Allard \cite{Allard}, Michael and Simon \cite{Michael-Simon}, and Castillon \cite{Castillon}. In particular, the Michael-Simon-Sobolev inequality implies an isoperimetric inequality for minimal surfaces. We refer to \cite{Brendle2} for a recent survey on geometric inequalities for minimal surfaces.

Finally, we refer to \cite{Cordero-Erausquin-McCann-Schmuckenschlager}, \cite{Gangbo-McCann}, \cite{McCann} for some of the earlier work on optimal transport and its applications to geometric inequalities.

\section{Proof of Theorem \ref{main.theorem}}

Let $\Sigma$ be a compact $n$-dimensional submanifold of $\mathbb{R}^{n+m}$, possibly with boundary $\partial \Sigma$. We denote by $g$ the Riemannian metric on $\Sigma$ and by $d(\cdot,\cdot)$ the Riemannian distance. For each point $x \in \Sigma$, we denote by $I\!I(x): T_x \Sigma \times T_x \Sigma \to T_x^\perp \Sigma$ the second fundamental form of $\Sigma$. As usual, the mean curvature vector $H(x) \in T_x^\perp \Sigma$ is defined as the trace of the second fundamental form.

We first consider the special case when $|\Sigma| = 1$. Let $\mu$ denote the Riemannian measure on $\Sigma$. We define a Borel measure $\nu$ on the unit ball $\bar{B}^{n+m}$ by 
\[\nu(G) = \int_G \rho(|\xi|^2) \, d\xi\] 
for every Borel set $G \subset \bar{B}^{n+m}$. With this understood, $\mu$ is a probability measure on $\Sigma$ and $\nu$ is a probability measure on $\bar{B}^{n+m}$. Let $\mathcal{J}$ denote the set of all pairs $(u,h)$ such that $u$ is an integrable function on $\Sigma$, $h$ is an integrable function on $\bar{B}^{n+m}$, and 
\begin{equation} 
\label{class.J}
u(x)-h(\xi)-\langle x,\xi \rangle \geq 0 
\end{equation}
for all $x \in \Sigma$ and all $\xi \in \bar{B}^{n+m}$. By Theorem 5.10 (iii) in \cite{Villani}, we can find a pair $(u,h) \in \mathcal{J}$ which maximizes the functional 
\begin{equation} 
\label{functional}
\int_{\bar{B}^{n+m}} h \, d\nu - \int_\Sigma u \, d\mu. 
\end{equation}
In fact, the result in \cite{Villani} shows that the maximizer $(u,h)$ may be chosen in such a way that $h$ is Lipschitz continuous and 
\begin{equation}
\label{c.transform} 
u(x) = \sup_{\xi \in \bar{B}^{n+m}} (h(\xi) + \langle x,\xi \rangle) 
\end{equation}
for all $x \in \Sigma$. 

Note that our notation differs from the one in \cite{Villani}. In our setting, the space $X$ is the unit ball $\bar{B}^{n+m}$ equipped with the measure $\nu$; the space $Y$ is the submanifold $\Sigma$ equipped with the Riemannian measure $\mu$; the cost function is given by $c(x,\xi) = -\langle x,\xi \rangle$ for $x \in \Sigma$ and $\xi \in \bar{B}^{n+m}$; the function $\psi$ in \cite{Villani} corresponds to the function $-h$; and the function $\phi$ in \cite{Villani} corresponds to the function $-u$ in this paper. The fact that $\psi$ can be chosen to be a $c$-convex function implies that $h$ is Lipschitz continuous (see \cite{Villani}, Definition 5.2). The fact that $\phi$ can be taken as the $c$-transform of $\psi$ corresponds to the statement (\ref{c.transform}) above (see \cite{Villani}, Definition 5.2).

It follows from (\ref{c.transform}) that $u$ is the restriction to $\Sigma$ of a convex function on $\mathbb{R}^{n+m}$ which is Lipschitz continuous with Lipschitz constant at most $1$. In particular, $u$ is Lipschitz continuous with Lipschitz constant at most $1$. Moreover, $u$ is semiconvex with a quadratic modulus of semiconvexity (see \cite{Villani}, Definition 10.10 and Example 10.11).

\begin{lemma}
\label{measure}
Let $E$ be a compact subset of $\Sigma$. Moreover, suppose that $G$ is a compact subset of $\bar{B}^{n+m}$ such that $u(x) - h(\xi) - \langle x,\xi \rangle > 0$ for all $x \in E$ and all $\xi \in \bar{B}^{n+m} \setminus G$. Then $\mu(E) \leq \nu(G)$.
\end{lemma} 

\textbf{Proof.} 
For every positive integer $j$, we define a compact set $G_j \subset \bar{B}^{n+m}$ by 
\[G_j = \{\xi \in \bar{B}^{n+m}: \text{\rm $\exists \, x \in E$ with $u(x) - h(\xi) - \langle x,\xi \rangle \leq j^{-1}$}\}.\] 
We define an integrable function $u_j$ on $\Sigma$ by $u_j = u - j^{-1} \cdot 1_E$. Moreover, we define an integrable function $h_j$ on $\bar{B}^{n+m}$ by $h_j = h - j^{-1} \cdot 1_{G_j}$. Using (\ref{class.J}), it is straightforward to verify that 
\[u_j(x)-h_j(\xi)-\langle x,\xi \rangle \geq 0\] 
for all $x \in \Sigma$ and all $\xi \in \bar{B}^{n+m}$. Therefore, $(u_j,h_j) \in \mathcal{J}$ for each $j$. Since the pair $(u,h)$ maximizes the functional (\ref{functional}), we obtain 
\[\int_{\bar{B}^{n+m}} h_j \, d\nu - \int_\Sigma u_j \, d\mu \leq \int_{\bar{B}^{n+m}} h \, d\nu - \int_\Sigma u \, d\mu\] 
for each $j$. This implies $\mu(E) \leq \nu(G_j)$ for each $j$. 

Finally, we pass to the limit as $j \to \infty$. Note that $G_{j+1} \subset G_j$ for each $j$. Since $E$ is compact and $u$ is continuous, we obtain 
\[\bigcap_{j=1}^\infty G_j \subset \{\xi \in \bar{B}^{n+m}: \text{\rm $\exists \, x \in E$ with $u(x) - h(\xi) - \langle x,\xi \rangle \leq 0$}\} \subset G.\] 
Putting these facts together, we conclude that 
\[\mu(E) \leq \lim_{j \to \infty} \nu(G_j) \leq \nu(G).\] 
This completes the proof of Lemma \ref{measure}. \\

Let us fix a large positive constant $K$ such that $|\langle x-\bar{x},y \rangle| \leq K \, d(x,\bar{x})^2$ for all points $x,\bar{x} \in \Sigma$ and all $y \in T_{\bar{x}}^\perp \Sigma$ with $|y| \leq 1$. For each point $\bar{x} \in \Sigma$, we define
\[\partial u(\bar{x}) = \{z \in T_{\bar{x}} \Sigma: \text{$u(x) - u(\bar{x}) - \langle x-\bar{x},z \rangle \geq -K \, d(x,\bar{x})^2$ for all $x \in \Sigma$}\}.\] 
We refer to $\partial u(\bar{x})$ as the subdifferential of $u$ at the point $\bar{x}$. 

\begin{lemma}
\label{subdifferential}
Fix a point $\bar{x} \in \Sigma$ and let $\xi \in \bar{B}^{n+m}$. Let $\xi^{\text{\rm tan}}$ denote the orthogonal projection of $\xi$ to the tangent space $T_{\bar{x}} \Sigma$. If $u(\bar{x}) - h(\xi) - \langle \bar{x},\xi \rangle = 0$, then $\xi^{\text{\rm tan}} \in \partial u(\bar{x})$. 
\end{lemma} 

\textbf{Proof.} 
By assumption, 
\[u(\bar{x}) - h(\xi) - \langle \bar{x},\xi \rangle = 0.\] 
Since 
\[u(x) - h(\xi) - \langle x,\xi \rangle \geq 0\] 
for all $x \in \Sigma$, it follows that 
\begin{equation} 
\label{ineq.1}
u(x) - u(\bar{x}) - \langle x-\bar{x},\xi \rangle \geq 0 
\end{equation}
for all $x \in \Sigma$. Using the fact that $\xi - \xi^{\text{\rm tan}} \in T_{\bar{x}}^\perp \Sigma$ and $|\xi - \xi^{\text{\rm tan}}| \leq |\xi| \leq 1$, we obtain 
\begin{equation} 
\label{ineq.2}
\langle x-\bar{x},\xi-\xi^{\text{\rm tan}} \rangle \geq -K \, d(x,\bar{x})^2 
\end{equation}
by our choice of $K$. Combining (\ref{ineq.1}) and (\ref{ineq.2}), we conclude that  
\begin{equation} 
u(x) - u(\bar{x}) - \langle x-\bar{x},\xi^{\text{\rm tan}} \rangle \geq -K \, d(x,\bar{x})^2. 
\end{equation}
Therefore, $\xi^{\text{\rm tan}} \in \partial u(\bar{x})$. This completes the proof of Lemma \ref{subdifferential}. \\

By Rademacher's theorem, $u$ is differentiable almost everywhere. At each point where $u$ is differentiable, the norm of its gradient is at most $1$. By Alexandrov's theorem (see Theorem 14.1 and Theorem 14.25 in \cite{Villani}), $u$ admits a Hessian in the sense of Alexandrov at almost every point. 

In the following, we fix a point $\bar{x} \in \Sigma \setminus \partial \Sigma$ with the property that $u$ admits a Hessian in the sense of Alexandrov at $\bar{x}$. Let $\hat{u}$ be a smooth function on $\Sigma$ such that $|u(x) - \hat{u}(x)| \leq o(d(x,\bar{x})^2)$ as $x \to \bar{x}$. 

Let us fix a small positive real number $\bar{r}$ so that $\frac{\sqrt{n}}{2} \, \bar{r} < d(\bar{x},\partial \Sigma)$ and $\frac{\sqrt{n}}{2} \, \bar{r}$ is smaller than the injectivity radius at $\bar{x}$. 

For each $r \in (0,\bar{r})$, we denote by $\hat{\omega}(r)$ the smallest nonnegative real number $\omega$ with the property that $|z - \nabla^\Sigma \hat{u}(x)| \leq \omega$ whenever $x \in \Sigma$, $z \in \partial u(x)$, and $d(x,\bar{x}) \leq \frac{\sqrt{n}}{2} \, r$.

\begin{lemma}
\label{properties.of.hat.omega}
The function $\hat{\omega}: (0,\bar{r}) \to [0,\infty)$ is monotone increasing and $\lim_{r \to 0} \frac{\hat{\omega}(r)}{r} = 0$. 
\end{lemma} 

\textbf{Proof.} 
The first statement follows immediately from the definition. The second property follows from the basic properties of the Alexandrov Hessian; see \cite{Villani}, Theorem 14.25 (i'). This completes the proof of Lemma \ref{properties.of.hat.omega}. \\

For each $r \in (0,\bar{r})$, we denote by $\hat{\delta}(r)$ the smallest nonnegative real number $\delta$ with the property that $D_\Sigma^2 \hat{u}(x) - \langle I\!I(x),\xi \rangle \geq -\delta \, g$ whenever $x \in \Sigma$, $\xi \in \bar{B}^{n+m}$, $u(x) - h(\xi) - \langle x,\xi \rangle = 0$, and $d(x,\bar{x}) \leq \frac{\sqrt{n}}{2} \, r$. 

\begin{lemma}
\label{properties.of.hat.delta}
The function $\hat{\delta}: (0,\bar{r}) \to [0,\infty)$ is monotone increasing and $\lim_{r \to 0} \hat{\delta}(r) = 0$. 
\end{lemma} 

\textbf{Proof.} 
The first statement follows immediately from the definition. To prove the second statement, we argue by contradiction. Suppose that $\limsup_{r \to 0} \hat{\delta}(r) > 0$. Then we can find a positive real number $\delta_0$, a sequence of points $x_j \in \Sigma$, and a sequence $\xi_j \in \bar{B}^{n+m}$ with the following properties: 
\begin{itemize}
\item $x_j \to \bar{x}$.
\item $u(x_j) - h(\xi_j) - \langle x_j,\xi_j \rangle = 0$ for each $j$. 
\item For each $j$, the first eigenvalue of $D_\Sigma^2 \hat{u}(x_j) - \langle I\!I(x_j),\xi_j \rangle$ is less than $-\delta_0$. 
\end{itemize}
After passing to a subsequence, we may assume that the sequence $\xi_j$ converges to $\bar{\xi} \in \bar{B}^{n+m}$. Since $\hat{u}$ is a smooth function, it follows that the first eigenvalue of $D_\Sigma^2 \hat{u}(\bar{x}) - \langle I\!I(\bar{x}),\bar{\xi} \rangle$ is strictly negative. Moreover, 
\[u(\bar{x}) - h(\bar{\xi}) - \langle \bar{x},\bar{\xi} \rangle = 0.\] 
Since 
\[u(x) - h(\bar{\xi}) - \langle x,\bar{\xi} \rangle \geq 0\] 
for all $x \in \Sigma$, it follows that 
\[u(x) - u(\bar{x}) - \langle x-\bar{x},\bar{\xi} \rangle \geq 0\] 
for all $x \in \Sigma$. Since $|u(x) - \hat{u}(x)| \leq o(d(x,\bar{x})^2)$ as $x \to \bar{x}$, we conclude that 
\[\hat{u}(x) - \hat{u}(\bar{x}) - \langle x-\bar{x},\bar{\xi} \rangle \geq -o(d(x,\bar{x}))^2\] 
as $x \to \bar{x}$. This implies $D_\Sigma^2 \hat{u}(\bar{x}) - \langle I\!I(\bar{x}),\bar{\xi} \rangle \geq 0$. This is a contradiction. This completes the proof of Lemma \ref{properties.of.hat.delta}. \\

Let $\{e_1,\hdots,e_n\}$ be an orthonormal basis of $T_{\bar{x}} \Sigma$. For each $r \in (0,\bar{r})$, we consider the cube 
\[W_r = \Big \{ z \in T_{\bar{x}} \Sigma: \max_{1 \leq i \leq n} |\langle z,e_i \rangle| \leq \frac{1}{2} \, r \Big \}.\] 
We denote by 
\[E_r = \exp_{\bar{x}}(W_r) \subset \Big \{ x \in \Sigma: d(x,\bar{x}) \leq \frac{\sqrt{n}}{2} \, r \Big \}\] 
the image of the cube $W_r$ under the exponential map. We further define 
\begin{align*} 
A_r 
&= \{(x,y): x \in E_r, \, y \in T_x^\perp \Sigma, \, |\nabla^\Sigma \hat{u}(x)|^2+|y|^2 \leq (1+\hat{\omega}(r))^2, \\ 
&\hspace{55mm} D_\Sigma^2 \hat{u}(x) - \langle I\!I(x),y \rangle \geq -\hat{\delta}(r) \, g\}. 
\end{align*} 
Clearly, $E_r$ is a compact subset of $\Sigma$ and $A_r$ is a compact subset of the normal bundle of $\Sigma$. We define a smooth map $\Phi: T^\perp \Sigma \to \mathbb{R}^{n+m}$ by 
\[\Phi(x,y) = \nabla^\Sigma \hat{u}(x) + y\] 
for $x \in \Sigma$ and $y \in T_x^\perp \Sigma$. Moreover, we denote by 
\[G_r = \{\xi \in \bar{B}^{n+m}: \text{\rm $\exists \, (x,y) \in A_r$ with $|\xi - \Phi(x,y)| \leq \hat{\omega}(r)$}\}\] 
the intersection of $\bar{B}^{n+m}$ with the tubular neighborhood of $\Phi(A_r)$ of radius $\hat{\omega}(r)$. Clearly, $G_r$ is a compact subset of $\bar{B}^{n+m}$.

\begin{lemma} 
\label{aux.1} 
Let $r \in (0,\bar{r})$. Then $u(x) - h(\xi) - \langle x,\xi \rangle > 0$ for all $x \in E_r$ and all $\xi \in \bar{B}^{n+m} \setminus G_r$. 
\end{lemma} 

\textbf{Proof.} 
We argue by contradiction. Suppose that there is a point $x \in E_r$ and a point $\xi \in \bar{B}^{n+m} \setminus G_r$ such that $u(x) - h(\xi) - \langle x,\xi \rangle = 0$. Let $\xi^{\text{\rm tan}}$ denote the orthogonal projection of $\xi$ to the tangent space $T_x \Sigma$. By Lemma \ref{subdifferential}, $\xi^{\text{\rm tan}} \in \partial u(x)$. Since $d(x,\bar{x}) \leq \frac{\sqrt{n}}{2} \, r$, it follows that  
\[|\xi^{\text{\rm tan}} - \nabla^\Sigma \hat{u}(x)| \leq \hat{\omega}(r)\] 
by definition of $\hat{\omega}(r)$. Let $y = \xi - \xi^{\text{\rm tan}} \in T_x^\perp \Sigma$. Then 
\[|\xi - \Phi(x,y)| = |\xi - \nabla^\Sigma \hat{u}(x) - y| = |\xi^{\text{\rm tan}} - \nabla^\Sigma \hat{u}(x)| \leq \hat{\omega}(r).\] 
Using the triangle inequality, we obtain  
\[\sqrt{|\nabla^\Sigma \hat{u}(x)|^2+|y|^2} = |\Phi(x,y)| \leq |\xi| + \hat{\omega}(r) \leq 1 + \hat{\omega}(r).\] 
Finally, since $d(x,\bar{x}) \leq \frac{\sqrt{n}}{2} \, r$, it follows that 
\[D_\Sigma^2 \hat{u}(x) - \langle I\!I(x),y \rangle = D_\Sigma^2 \hat{u}(x) - \langle I\!I(x),\xi \rangle \geq -\hat{\delta}(r) \, g\] 
by definition of $\hat{\delta}(r)$. To summarize, we have shown that $(x,y) \in A_r$ and $|\xi - \Phi(x,y)| \leq \hat{\omega}(r)$. Consequently, $\xi \in G_r$, contrary to our assumption. This completes the proof of Lemma \ref{aux.1}. \\

\begin{lemma}
\label{aux.2}
Let $r \in (0,\bar{r})$. Then $\mu(E_r) \leq \nu(G_r)$. 
\end{lemma} 

\textbf{Proof.} 
This follows by combining Lemma \ref{measure} and Lemma \ref{aux.1}. \\

\begin{proposition} 
\label{Monge.Ampere}
Fix a point $\bar{x} \in \Sigma \setminus \partial \Sigma$ with the property that $u$ admits a Hessian in the sense of Alexandrov at $\bar{x}$. Let $\hat{u}$ be a smooth function on $\Sigma$ such that $|u(x) - \hat{u}(x)| \leq o(d(x,\bar{x})^2)$ as $x \to \bar{x}$. Let 
\[S = \{y \in T_{\bar{x}}^\perp \Sigma: |\nabla^\Sigma \hat{u}(\bar{x})|^2+|y|^2 \leq 1, \, D_\Sigma^2 \hat{u}(\bar{x}) - \langle I\!I(\bar{x}),y \rangle \geq 0\}.\] 
Then 
\[1 \leq \int_S \det (D_\Sigma^2 \hat{u}(\bar{x}) - \langle I\!I(\bar{x}),y \rangle) \, \rho(|\nabla^\Sigma \hat{u}(\bar{x})|^2+|y|^2) \, dy.\] 
\end{proposition}

\textbf{Proof.} 
In the following, we fix an arbitrary positive integer $j$. We define 
\[S_j = \{y \in T_{\bar{x}}^\perp \Sigma: |\nabla^\Sigma \hat{u}(\bar{x})|^2+|y|^2 \leq 1+j^{-1}, \, D_\Sigma^2 \hat{u}(\bar{x}) - \langle I\!I(\bar{x}),y \rangle \geq -j^{-1} \, g\}.\] 
For each $r \in (0,\bar{r})$, we decompose the normal space $T_{\bar{x}}^\perp \Sigma$ into compact cubes of size $r$. Let $\mathcal{Q}_r$ denote the collection of all the cubes in this decomposition. Moreover, we denote by $\mathcal{Q}_{r,j} \subset \mathcal{Q}_r$ the set of all cubes in $\mathcal{Q}_r$ that are contained in the set $S_j$. We define a smooth map 
\[\Psi: W_r \times T_{\bar{x}}^\perp \Sigma \to \mathbb{R}^{n+m}, \, (z,y) \mapsto \Phi(\exp_{\bar{x}}(z),P_z y),\] 
where $P_z: T_{\bar{x}}^\perp \Sigma \to T_{\exp_{\bar{x}}(z)}^\perp \Sigma$ denotes the parallel transport along the geodesic $t \mapsto \exp_{\bar{x}}(tz)$ (see \cite{ONeill}, pp.~114--115). Since $\lim_{r \to 0} \hat{\omega}(r) = 0$ and $\lim_{r \to 0} \hat{\delta}(r) = 0$, we obtain 
\[\Phi(A_r) \subset \bigcup_{Q \in \mathcal{Q}_{r,j}} \Psi(W_r \times Q),\] 
provided that $r$ is sufficiently small (depending on $j$). This implies 
\begin{align*} 
G_r 
&= \{\xi \in \bar{B}^{n+m}: \text{\rm $\exists \, (x,y) \in A_r$ with $|\xi - \Phi(x,y)| \leq \hat{\omega}(r)$}\} \\ 
&\subset \bigcup_{Q \in \mathcal{Q}_{r,j}} \{\xi \in \bar{B}^{n+m}: \text{\rm $\exists \, (z,y) \in W_r \times Q$ with $|\xi - \Psi(z,y)| \leq \hat{\omega}(r)$}\}, 
\end{align*}
provided that $r$ is sufficiently small (depending on $j$).

We next observe that 
\[|\det D\Psi(0,y)| = |\det D\Phi(\bar{x},y)| = |\det(D_\Sigma^2 \hat{u}(\bar{x}) - \langle I\!I(\bar{x}),y \rangle)|\] 
for all $y \in T_{\bar{x}}^\perp \Sigma$. Hence, if $r$ is sufficiently small (depending on $j$), then we obtain 
\begin{align} 
\label{change.of.variables.formula}
&\nu \big ( \{\xi \in \bar{B}^{n+m}: \text{\rm $\exists \, (z,y) \in W_r \times Q$ with $|\xi - \Psi(z,y)| \leq \hat{\omega}(r)$}\} \big ) \notag \\ 
&\leq r^n \int_Q \big [ |\det(D_\Sigma^2 \hat{u}(\bar{x}) - \langle I\!I(\bar{x}),y \rangle)| \, \rho(|\nabla^\Sigma \hat{u}(\bar{x})|^2+|y|^2) + j^{-1} \big ] \, dy 
\end{align} 
for each cube $Q \in \mathcal{Q}_{r,j}$. To justify (\ref{change.of.variables.formula}), we argue as in the proof of the classical change-of-variables formula  (see \cite{Rudin}, pp.~150--156). We also use the fact that $\lim_{r \to 0} \frac{\hat{\omega}(r)}{r} = 0$. 

Summation over all cubes $Q \in \mathcal{Q}_{r,j}$ gives 
\begin{align*} 
&\nu(G_r) \\ 
&\leq \sum_{Q \in \mathcal{Q}_{r,j}} \nu \big ( \{\xi \in \bar{B}^{n+m}: \text{\rm $\exists \, (z,y) \in W_r \times Q$ with $|\xi - \Psi(z,y)| \leq \hat{\omega}(r)$}\} \big ) \\ 
&\leq r^n \int_{S_j} \big [ |\det(D_\Sigma^2 \hat{u}(\bar{x}) - \langle I\!I(\bar{x}),y \rangle)| \, \rho(|\nabla^\Sigma \hat{u}(\bar{x})|^2+|y|^2) + j^{-1} \big ] \, dy, 
\end{align*}
provided that $r$ is sufficiently small (depending on $j$).

On the other hand, Lemma \ref{aux.2} implies that $\mu(E_r) \leq \nu(G_r)$ for each $r \in (0,\bar{r})$. Thus, we conclude that 
\begin{align*} 
1 &= \limsup_{r \to 0} r^{-n} \, \mu(E_r) \\ 
&\leq \limsup_{r \to 0} r^{-n} \, \nu(G_r) \\ 
&\leq \int_{S_j} \big [ |\det(D_\Sigma^2 \hat{u}(\bar{x}) - \langle I\!I(\bar{x}),y \rangle)| \, \rho(|\nabla^\Sigma \hat{u}(\bar{x})|^2+|y|^2) + j^{-1} \big ] \, dy. 
\end{align*} 
Finally, we pass to the limit as $j \to \infty$. Note that $S_{j+1} \subset S_j$ for each $j$. Moreover, $\bigcap_{j=1}^\infty S_j = S$. This gives 
\[1 \leq \int_S |\det (D_\Sigma^2 \hat{u}(\bar{x}) - \langle I\!I(\bar{x}),y \rangle)| \, \rho(|\nabla^\Sigma \hat{u}(\bar{x})|^2+|y|^2) \, dy.\] 
Since $D_\Sigma^2 \hat{u}(\bar{x}) - \langle I\!I(\bar{x}),y \rangle \geq 0$ for all $y \in S$, the assertion follows. This completes the proof of Proposition \ref{Monge.Ampere}. \\

\begin{corollary} 
\label{estimate.for.Laplacian}
Fix a point $\bar{x} \in \Sigma \setminus \partial \Sigma$ with the property that $u$ admits a Hessian in the sense of Alexandrov at $\bar{x}$. Let $\hat{u}$ be a smooth function on $\Sigma$ such that $|u(x) - \hat{u}(x)| \leq o(d(x,\bar{x})^2)$ as $x \to \bar{x}$. Then 
\[n \, \alpha^{-\frac{1}{n}} \leq \Delta_\Sigma \hat{u}(\bar{x}) + |H(\bar{x})|,\] 
where $\alpha$ is defined by (\ref{def.alpha}). 
\end{corollary}

\textbf{Proof.} 
We argue by contradiction. If the assertion is false, then there exists a real number $\hat{\alpha} > \alpha$ such that 
\[\Delta_\Sigma \hat{u}(\bar{x}) + |H(\bar{x})| \leq n \, \hat{\alpha}^{-\frac{1}{n}}.\] 
Let 
\[S = \{y \in T_{\bar{x}}^\perp \Sigma: |\nabla^\Sigma \hat{u}(\bar{x})|^2+|y|^2 \leq 1, \, D_\Sigma^2 \hat{u}(\bar{x}) - \langle I\!I(\bar{x}),y \rangle \geq 0\}.\] 
The arithmetic-geometric mean inequality gives 
\[0 \leq \det(D_\Sigma^2 \hat{u}(\bar{x}) - \langle I\!I(\bar{x}),y \rangle) \leq \Big ( \frac{\Delta_\Sigma \hat{u}(\bar{x}) - \langle H(\bar{x}),y \rangle}{n} \Big )^n \leq \hat{\alpha}^{-1}\] 
for all $y \in S$. Using Proposition \ref{Monge.Ampere}, we obtain 
\begin{align*} 
1 
&\leq \int_S \det(D_\Sigma^2 \hat{u}(\bar{x}) - \langle I\!I(\bar{x}),y \rangle) \, \rho(|\nabla^\Sigma \hat{u}(\bar{x})|^2+|y|^2) \, dy \\ 
&\leq \int_S \hat{\alpha}^{-1} \, \rho(|\nabla^\Sigma \hat{u}(\bar{x})|^2+|y|^2) \, dy \\ 
&\leq \hat{\alpha}^{-1} \, \alpha. 
\end{align*} 
In the last step, we have used the definition of $\alpha$; see (\ref{def.alpha}). Thus $\hat{\alpha} \leq \alpha$, contrary to our assumption. This completes the proof of Corollary \ref{estimate.for.Laplacian}. \\

After these preparations, we may now complete the proof of Theorem \ref{main.theorem}. Corollary \ref{estimate.for.Laplacian} implies that 
\begin{equation} 
\label{pointwise.estimate}
n \, \alpha^{-\frac{1}{n}} \leq \Delta_\Sigma u + |H| 
\end{equation} 
almost everywhere, where $\Delta_\Sigma u$ denotes the trace of the Alexandrov Hessian of $u$. The distributional Laplacian of $u$ may be decomposed into its singular and absolutely continuous part. By Alexandrov's theorem (see Theorem 14.1 in \cite{Villani}), the density of the absolutely continuous part is given by the trace of the Alexandrov Hessian of $u$. The singular part of the distributional Laplacian of $u$ is nonnegative since $u$ is semiconvex. This implies 
\begin{equation} 
\label{integration.by.parts}
\int_\Sigma \eta \, \Delta_\Sigma u \leq -\int_\Sigma \langle \nabla^\Sigma \eta,\nabla^\Sigma u \rangle 
\end{equation} 
for every nonnegative smooth function $\eta: \Sigma \to \mathbb{R}$ that vanishes in a neighborhood of $\partial \Sigma$. Combining (\ref{pointwise.estimate}) and (\ref{integration.by.parts}), we obtain 
\begin{align*} 
n \, \alpha^{-\frac{1}{n}} \int_\Sigma \eta 
&\leq \int_\Sigma \eta \, \Delta_\Sigma u + \int_\Sigma \eta \, |H| \\ 
&\leq -\int_\Sigma \langle \nabla^\Sigma \eta,\nabla^\Sigma u \rangle + \int_\Sigma \eta \, |H| \\ 
&\leq \int_\Sigma |\nabla^\Sigma \eta| + \int_\Sigma \eta \, |H| 
\end{align*} 
for every nonnegative smooth function $\eta: \Sigma \to \mathbb{R}$ that vanishes in a neighborhood of $\partial \Sigma$. By a straightforward limiting procedure, this implies 
\[n \, \alpha^{-\frac{1}{n}} \, |\Sigma| \leq |\partial \Sigma| + \int_\Sigma |H|.\] 
This completes the proof of Theorem \ref{main.theorem} in the special case when $|\Sigma| = 1$. The general case follows by scaling.

\section{Proof of Corollary \ref{consequence.of.main.theorem}}

In this final section, we explain how Corollary \ref{consequence.of.main.theorem} follows from Theorem \ref{main.theorem}. Assume that $n \geq 2$ and $m \geq 2$. We can find a find a sequence of continuous functions $\rho_j: [0,\infty) \to (0,\infty)$ such that $\int_{\bar{B}^{n+m}} \rho_j(|\xi|^2) \, d\xi = 1$, 
\[\sup_{[0,1-j^{-1}]} \rho_j \leq o(1),\] 
and 
\[\sup_{[1-j^{-1},1]} \rho_j \leq \frac{2j}{(n+m) \, |B^{n+m}|} + o(j)\] 
as $j \to \infty$. For each point $z \in \mathbb{R}^n$, we obtain 
\begin{align*} 
&\int_{\{y \in \mathbb{R}^m: \, |z|^2+|y|^2 \leq 1\}} \rho_j(|z|^2+|y|^2) \, dy \\ 
&\leq |B^m| \, (1-|z|^2-j^{-1})_+^{\frac{m}{2}} \, \sup_{[0,1-j^{-1}]} \rho_j \\ 
&+ |B^m| \, \big [ (1-|z|^2)_+^{\frac{m}{2}} - (1-|z|^2-j^{-1})_+^{\frac{m}{2}} \big ] \, \sup_{[1-j^{-1},1]} \rho_j \\ 
&\leq |B^m| \, \sup_{[0,1-j^{-1}]} \rho_j + \frac{m}{2} \, |B^m| \, j^{-1} \, \sup_{[1-j^{-1},1]} \rho_j.
\end{align*} 
In the last step, we have used the fact that $m \geq 2$. This implies 
\[\sup_{z \in \mathbb{R}^n} \int_{\{y \in \mathbb{R}^m: \, |z|^2+|y|^2 \leq 1\}} \rho_j(|z|^2+|y|^2) \, dy \leq \frac{m \, |B^m|}{(n+m) \, |B^{n+m}|} + o(1)\] 
as $j \to \infty$. Therefore, the assertion follows from Theorem \ref{main.theorem}.

\end{document}